\documentclass[a4paper,12pt]{article}
\usepackage[LCY,T2A]{fontenc}
\usepackage[cp1251]{inputenc}
\usepackage{amsthm}
\usepackage{amsmath}
\usepackage{amssymb}
\usepackage{mathrsfs}
\usepackage[english,russian]{babel}
\usepackage[active]{srcltx}

\begin{document}

\begin{center}
{\bf One Counterexample of Comonotone Approximation \\ of
$2\pi$-periodic Function on Trigonometric Polynomials}
\end{center}
\author{M.~G.~Pleshakov}%

\begin{center}
Pleshakov Mikhail Gennad'evich, Saratov State University,
Russia,\\
410012, Saratov, Astrakhanskaya st., 83, pleshakovmg@gmail.com
\end{center}

{\bf Abstract.} Let $2s$ points $y_i=-\pi\le y_{2s}<\ldots<y_1<\pi$
be given. Using these points, we define the points $y_i$ for all
integer indices $i$ by the equality $y_i=y_{i+2s}+2\pi$. We shall
write $f\in\bigtriangleup^{(1)}(Y)$ if $f$ is a $2\pi$-periodic
function and $f$ does not decrease on $[y_i, y_{i-1}]$ if $i$ is
odd; and $f$ does not increase on $[y_i, y_{i-1}]$ if $i$ is even.
We denote $E_n^{(1)}(f;Y)$ the value of the best uniform
comonotone approximation.  In this article the following
counterexample of comonotone approximation is proved.\newline
 \textbf{Example.} {\it For each $k\in\Bbb
N$, $k>3$, and $n\in\Bbb N$ there a function $f(x):=f(x;s,Y,n,k)$
exists, such that $f\in\bigtriangleup^{(1)}(Y)\bigcap{\Bbb
C}^{(1)}$ and
$$ E_n^{(1)}(f;Y)>B_Yn^{\frac k3 -1}\frac 1n\omega_k\left(f';\frac 1n\right),
                                                                  $$}
where $B_Y=$const, depending only on $Y$ and $k$; $\omega_k$ is
the modulus of smoothness of order $k$, of $f$.%


\smallskip
Получение оценки уклонения при равномерном приближении непрерывных
функций алгебраическими многочленами и
тригоно\-ме\-т\-ри\-че\-с\-ки\-ми полиномами
 является одной из основных задач в теории приближения функций.
Наиболее широкое применение в теоретических исследованиях и в
прикладных областях математики получили неравенства типа
Джексона-Зигмунда-Стечкина [1-3], типа
Никольского-Ти\-м\-а\-на-Дзя\-ды\-ка-Фройда-Теляковского-Брудного
[4-9]. Особый интерес представляет случай, когда приближение
является формосохраняющим (Shape-preserving Approximation), т.е.
когда аппарат приближения сохраняет некоторые свойства
приближаемой функции (монотонность, выпуклость и т.д.) В 1969 году
G.G. Lorentz и K.L. Zeller [10] построили пример, который
показывает, что величина наилучшего монотонного приближения
алгебраическими многочленами монотонной функции по порядку вообще
говоря ``хуже'' величины наилучшего приближения без ограничений. В
1979 году А.С. Шведов [11, 12] построил пример, показывающий, что
оценка типа Джексона-Стечкина величины приближения монотонной
функции монотонными многочленами через модуль непрерывности
порядка 3 и выше вообще неверна, в отличие от приближения без
ограничений.

Однако результаты по комонотонному приближению периодических
функций три\-г\-о\-н\-о\-ме\-т\-р\-и\-че\-с\-ки\-ми полиномами, за
исключением результата G.G. Lorentz и K.L. Zeller 1968 года,
касающегося так называемых ``колоколообразных'' функций долгое
время не были известны.

В данной статье построен контрпример, указывающий, что оценка типа
Джексона-Стечкина величины приближения кусочно-монотонной
периодической  функции комонотонными тригонометрическими
полиномами через модуль непрерывности порядка 3 и выше вообще
неверна.

Пусть $\Bbb C$ -- пространство непрерывных $2\pi$-периодических
действительнозначных функций $f$ с равномерной нормой $\Vert
f\Vert :=\max\limits_{x\in\Bbb R}\vert f(x)\vert $; $\omega(f;t) $
--  модуль непрерывности  функции  $f$; ${\Bbb T}_n,\ n\in\Bbb N$,
-- пространство тригонометрических полиномов
$$\tau_n(x):= a_0 + \sum\limits_{k=1}^n (a_k\cos kx + b_k\sin kx)$$
порядка $\leq n$.

Пусть на промежутке $[-\pi ,\pi )$ заданы $2s$ точек $y_i$:
$-\pi\leq y_{2s}<y_{2s-1}<\ldots <y_1<\pi$. Отправляясь от этих
точек, при помощи равенства $y_i=y_{i+2s}+2\pi$ определим точки
$y_i$ для всех целых индексов $i$; в частности, $y_0=y_{2s}+2\pi$,
$ y_{2s+1}=y_1-2\pi$ и т.д. Обозначим $Y:=\{ y_i\}_{i\in\Bbb Z}$.
Множество всех таких наборов обозначим ${\Bbb Y}_{2s}$. Будем
писать
$$f\in\bigtriangleup^{(1)}(Y),$$ если $f(x)$ -- $2\pi$-периодическая
непрерывная функция и $f(x)$ не убывает на $[y_i,y_{i-1}]$, если
$i$ нечетное; $f(x)$ не возрастает на $[y_i,y_{i-1}],$ если $i$
четное.

\medskip
Обозначим
$$\Pi(x):=\prod^{2s}_{i=1}\sin \frac {1}{2}(x-y_i),              $$
и заметим, что $\Pi\in{\Bbb T}_s,$ то есть $\Pi(x)$ --
тригонометрический полином порядка $s$.

Зафиксируем $s\in\Bbb N$ и набор $\{y_i\}_{i\in\Bbb Z}=Y\in{\Bbb
Y}_{2s}$. В силу периодичности без потери общности будем считать,
что точка $0$ принадлежит набору $Y$, т.е. $y_{i_*}=0$ при
некотором $i_*\in\Bbb Z$.

Обозначим
$$\Pi_*(x):=\prod_{i=1\ i\ne i_*}^{2s} \sin\frac {x-y_i}2.$$
Для определённости будем считать, что $i_*$ -- нечётное число.
Тогда
$$\Pi_*(0)>0.\eqno (1)$$

Обозначим через $2d$ расстояние от $y_{i_*}$ до ближайшей точки
набора $Y$, заметим,
$$
d\le\frac\pi 2,
$$
$$
\Pi_*(x)>0,\quad x\in(-2d,2d).\eqno (2)
$$
Положим
$$M:=\max_{x\in\Bbb R}|\Pi_*(x)|,\ \ \ \ M_1:=\max_{x\in\Bbb R}|\Pi_*'
(x)|,$$
$$m:=\min_{x\in [-d, d]}\Pi_*(x).$$

Отправляясь от набора $Y$, определим натуральное число $N$. А
имен\-но, обозначим через $N$ наименьшее из чисел, удовлетворяющих
неравенству
$$m\sin^3\frac  d8\ge\frac 5N\left(M+M_1\right).  \eqno (3) $$
Тогда
$$m\sin\frac  d8\ge\frac 5N\left(M+M_1\right),  \eqno (4) $$
следовательно
$$
d>\frac {40}N.                                  \eqno (5)
$$
Выберем натуральное число $j^*$ из условия
$$\frac {\pi}{N}+j^*\frac {2\pi}{N}\le d<\frac {\pi}{N}+(j^*+1)\frac
{2\pi}{N}.$$ Обозначим
$$d^*:=\frac {\pi}{N}+j^*\frac {2\pi}{N}$$
и заметим,
$$\frac 12d<d^*\le d.                            \eqno (6) $$

При построении контрпримера будет использовано ядро Джексона
$$J_N(t)=\frac 3{2N(2N^2+1)}\left(\frac {\sin\frac {Nt}2}{\sin\frac t2}
\right)^4. $$ Напомним (см. например [13, с. 127]),
не\-ко\-т\-о\-ры\-е сво\-й\-ст\-ва ядра Джексона:

\noindent а) $J_N(t)$ является чётным неотрицательным
тригонометрическим полиномом порядка $2N-2$;

\noindent б) $$\frac 1{\pi}\int\limits_{-\pi}^{\pi}J_N(t)dt=1;
\eqno (7)$$

\noindent в) для любой непрерывно дифференцируемой периодической
функции $g$ в каждой точке $x$ имеет место неравенство
$$\frac 1{\pi}\left\vert\int\limits_{-\pi}^{\pi}(g(t)-g(x))
J_N(t-x)dt\right\vert\le\frac 5N\|g'\|.                \eqno (8)$$
Обозначим
$$\tilde M:=\frac 1{\pi}\|J_N\|,
$$
$$
\tilde m:=\frac 1{\pi}\min_{t\in\left[-\frac {\pi}{2N}, \frac
{\pi}{2N} \right]}J_N(t-d^*)= \frac 1{\pi}\min_{t\in\left[-\frac
{\pi}{2N}, \frac {\pi}{2N} \right]}J_N(t+d^*),
$$
и заметим, что $\tilde m>0$. Наконец, положим
$$\overline M:=2+\pi^3\sqrt {\frac {M\tilde M}{m\tilde m}}.$$
Всюду далее в главе предполагаем, что число $b$ удовлетворяет
неравенствам
$$
0<b<\frac\pi{2N\overline M},                           \eqno (9)
$$
в частности, с учетом (5) и (6),
$$
\frac{d^*-2b}2>\frac d8.                             \eqno (10)
$$

Для доказательства примера нам потребуется ряд лемм.

\medskip
{\bf Л е м м а 2.} {\it Для любого $b$ найдётся  положительное
число $\alpha_b<1$ такое, что для функции

$$
Q(x):=Q(x;b)=\frac 1{\pi}\int\limits_{0}^{x}\sin\frac
{t-b}2\sin\frac {t+b}2\Pi(t)
\left(\alpha_bJ_N(t-d^*)+(1-\alpha_b)J_N(t+d^*)\right)dt
$$

имеет место равенство
            $$ Q(2\pi )=0.                                      \eqno (11) $$}
{\it Д о к а з а т е л ь с т в о.} Обозначим
$$Q_r(x,b):=Q_r(x):=\frac 1{\pi}\int\limits_{0}^{x}\sin\frac {t-b}2\sin\frac {t+b}2
\Pi(t)J_N(t-d^*)dt.$$ Представим $Q_r(2\pi )$ в виде

\begin{multline*} Q_r(2\pi):=\frac
1{\pi}\int\limits_{0}^{2\pi}\sin\frac {d^*-b}2\sin
\frac {d^*+b}2\Pi(d^*)J_N(t-d^*)dt+\\
+\frac 1{\pi}\int\limits_{0}^{2\pi}\left(\sin\frac {t-b}2\sin\frac
{t+b}2\Pi
(t)-\sin\frac {d^*-b}2\sin\frac {d^*+b}2\Pi(d^*)\right)J_N(t-d^*)dt\\
=:I_1+I_2.\end{multline*}

  В силу  (7) и (10)
$$
I_1=\sin\frac {d^*-b}2\sin\frac {d^*+b}2\Pi_*(d^*)\sin\frac
{d^*}2> m_0\sin^3\frac d8.
$$
Согласно (8)
$$|I_2|\le\frac 5N\left\Vert\left(\sin\frac {.-b}2\sin\frac {.+b}2\sin\frac
.2\Pi_*(\ .\ )\right)'\right\Vert\le\frac 5N\left(M+M_1\right).$$
Поэтому в силу (3)
$$Q_r(2\pi )=I_1+I_2>0.$$
Аналогично, обозначим
$$Q_l(x,b):=Q_l(x):=\frac 1{\pi}\int\limits_{0}^{x}\sin\frac {t-b}2\sin\frac {t+b}2
\Pi(t)J_N(t+d^*)dt.$$ Представим $Q_l(2\pi )$ в виде

\begin{multline*} Q_l(2\pi):=\frac
1{\pi}\int\limits_{0}^{2\pi}\sin\frac {-d^*-b}2\sin
\frac {-d^*+b}2\Pi(-d^*)J_N(t+d^*)dt+\\
+\frac 1{\pi}\int\limits_{0}^{2\pi}\left(\sin\frac {t-b}2\sin\frac
{t+b}2\Pi
(t)-\sin\frac {-d^*-b}2\sin\frac {-d^*+b}2\Pi(-d^*)\right)\cdot\\
\cdot J_N(t+d^*)dt=I_1+I_2.\end{multline*}

В силу  (7) и (10)
$$
I_1=\sin\frac {-d^*-b}2\sin\frac {-d^*+b}2\Pi^*(-d^*)\sin\frac
{-d^*}2< -m_0\sin^3\frac d8.
$$
Согласно (8)
$$|I_2|\le\frac 5N\left\Vert\left(\sin\frac {.-b}2\sin\frac {.+b}2\sin\frac
.2\Pi_*(\ .\ )\right)'\right\Vert\le\frac 5N\left(M+M_1\right).$$
Поэтому в силу (3)
$$Q_l(2\pi )=I_1+I_2<0.$$
Теперь осталось выбрать $\alpha_b$ из условия
$$\alpha_bQ_r(2\pi )+(1-\alpha_b)Q_l(2\pi )=0,$$
т.е.
$$\alpha_b=-\frac {Q_l(2\pi )}{Q_r(2\pi )-Q_l(2\pi )},$$
и заметить, что $0<\alpha_b< 1$.
 Лемма доказана.

\medskip
Равенство (11) означает, что $Q$ есть тригонометрический полином
порядка $2N+s-1$.

\medskip
{\it О п р е д е л е н и е 2.} При каждом $b$ назовём ``правым
$b$-корытом'' $2\pi$-периодическую функцию $K_{r,b}(x)$, имеющую
свойства:
$$
K_{r,b}\in {\Bbb C}^{(4)};
$$
$$
0\le K_{r,b}(x)\le 1,\quad x\in{\Bbb R};
$$
\[ K_{r, b}(x):=\left\{\begin{array}{rl} 0,\ \ \ \ \text {если}\ -\frac{\overline
Mb}2\le x\le b, \\
 1,\ \ \ \ \text {если } x\in[-\pi,-\overline Mb]\cup[2b,\pi]; \tag {12}
 \end{array}\right.\]

``левым $b$-корытом'' $2\pi$-периодическую функцию $K_{l,b}(x)$,
имеющую св\-о\-й\-с\-т\-ва:
$$
K_{l,b}\in {\Bbb C}^{(4)};
$$
$$
0\le K_{l,b}(x)\le 1,\quad x\in{\Bbb R};
$$
\[K_{l, b}(x):=\left\{\begin{array}{rl} 0,\ \ \ \ \mbox {если}\ -b\le x\le
\frac{\overline Mb}2,  \\ 1,\ \ \ \ \mbox {если}\
x\in[-\pi,-2b]\cup[\overline M,\pi].   \end{array}\right. \]

Замечание. В примере  достаточно, чтобы функция $K_{r,b}$ была
"просто" непрерывной. Поэтому  в лемме 2 в качестве $K_{r,b}$
можно взять, скажем, кусочно-линейную функцию. То же относится к
$K_{l,b}$.

\medskip
{\bf Л е м м а 3.} {\it Для любого $b$ имеет место неравенство} $$
\int\limits_{-\pi}^{\pi}K_{r,b}(t)\sin\frac {t-b}2\sin \frac
{t+b}2\Pi
(t)\left(\alpha_bJ_N(t-d^*)+(1-\alpha_b)J_N(t+d^*)\right)dt>0.
$$
{\it Д о к а з а т е л ь с т в о.} Вследствие леммы 2
$$
\int\limits_{-\pi}^{\pi}\sin\frac {t-b}2\sin\frac {t+b}2 \Pi
(t)\left(\alpha_bJ_N(t-d^*)+(1-\alpha_b)J_N(t+d^*)\right)dt=0,
$$
поэтому достаточно доказать неравенство

\begin{multline*}
\left\vert \int\limits_{-\frac{\overline Mb}2}^{-b}\sin\frac
{t-b}2\sin\frac {t+b}2
\Pi(t)\left(\alpha_bJ_N(t-d^*)+(1-\alpha_b)J_N(t+d^*)\right)dt\right\vert >\\
>\int\limits_{-b}^{2b}\sin\frac {t-b}2\sin\frac {t+b}2\Pi(t)
\left(\alpha_bJ_N(t-d^*)+(1-\alpha_b)J_N(t+d^*)\right)dt
.\end{multline*}

Если $-\frac{\overline Mb}2\le t\le -b$, то
$$\left\vert\sin\frac {t-b}2\right\vert \ge\sin b>\frac {2b}{\pi},$$
$$\left\vert\Pi_*(t)\right\vert \ge m,\ \ \
\left\vert\sin\frac t2\right\vert \ge\frac b{\pi},$$
$$\frac 1\pi J_N(t-d^*)\ge\tilde m,\ \ \ \frac 1\pi J_N(t+d^*)\ge\tilde m,$$
т.е.

\begin{multline*} \left\vert \int\limits_{-\frac{\overline
Mb}2}^{-b}\sin\frac {t-b}2\sin\frac {t+b}2
\Pi(t)\left(\alpha_bJ_N(t-d^*)+(1-\alpha_b)J_N(t+d^*)\right)dt\right\vert >\\
>\frac {2m\tilde mb^2}{\pi}\left\vert\int\limits_{-\frac{\overline Mb}2}^{-b}
\sin\frac {t+b}2dt\right\vert
=\frac {8m\tilde mb^2}{\pi}\sin^2\frac {(1-\frac{\overline M}2)b}4\ge\\
\ge\frac {2m\tilde m(1-\frac{\overline M}2)^2b^4}{\pi^3}.
\end{multline*}
С другой стороны
\begin{multline*} \int\limits_{-b}^{2b}\sin\frac
{t-b}2\sin\frac {t+b}2
\Pi(t)\left(\alpha_bJ_N(t-d^*)+(1-\alpha_b)J_N(t+d^*)\right)dt<\\
<3b(\sin b\sin\frac{3b}2\sin b)M\pi\tilde M<\frac{9\pi M\tilde
Mb^4}2.
\end{multline*}

В силу выбора $\overline M$ получаем
$$\frac {2m\tilde m(1-\frac{\overline M}2)^2b^4}{\pi^3}>
\frac {9\pi M\tilde Mb^4}2.$$ Лемма доказана.

\noindent Аналогично доказывается

{\bf Л е м м а 4.} {\it Для  любого $b$ имеет место неравенство
\begin{multline*} \int\limits_{-\pi}^{\pi}K_{l,b}(t)\sin\frac
{t-b}2\sin
\frac {t+b}2\Pi (t)\left(\alpha_bJ_N(t-d^*)+ \right .\\
\left .+(1-\alpha_b)J_N(t+d^*)\right)dt<0.
\end{multline*} }

Следствием предыдущих двух лемм является

{\bf Л е м м а 5.} {\it Для каждого $b$ существует
$2\pi$-периодическая $4$ разa непрерывно дифференцируемая функция
$K_b(x)$, которая имеет свойства: $K_b(x)=0$, если $|x|<b$;
$K_b(x)=1$, если $\overline Mb<|x|\le\pi$; $0\le K_b(x)\le 1,\quad
x\in{\Bbb R}$;
$$\int\limits_{0}^{2\pi}K_{b}(t)\sin\frac {t-b}2\sin
\frac
{t+b}2\Pi(t)\left(\alpha_bJ_N(t-d^*)+(1-\alpha_b)J_N(t+d^*)\right)dt=0.
                                                               \eqno (13) $$ }

{\it Д о к а з а т е л ь с т в о.} Из лемм  3 и 4 следует, что $$
I_1:=\frac 1{\pi}\int\limits_{-\pi}^{\pi}K_{r,b}(t)\sin\frac
{t-b}2 \sin\frac {t+b}2\Pi (t)\left(\alpha_bJ_N(t-d^*)+
(1-\alpha_b)J_N(t+d^*)\right)dt>0;  $$
$$I_2:=\frac
1{\pi}\int\limits_{-\pi}^{\pi}K_{l,b}(t)\sin\frac {t-b}2 \sin\frac
{t+b}2\Pi
(t)\left(\alpha_bJ_N(t-d^*)+(1-\alpha_b)J_N(t+d^*)\right)dt<0. $$
Остаётся выбрать $\gamma_b\in(0,1)$ из условия
$$
\gamma_bI_1+(1-\gamma_b)I_2=0
$$
и положить
$$
K_b(x):=\gamma_bK_{r,b}(x)+(1-\gamma_b)K_{l,b}(x).
$$
Лемма доказана.

Положим
$$q_b(x):=\sin\frac{(x-b)}2\sin\frac {(x+b)}2\sin\frac {x}2,$$
$$
g_b(x):= \frac 1{\pi}\int_{0}^xq_b(t)K_b(t)\frac {\Pi_*(t)}{M}
(\alpha_bJ_N(t-d^*)+(1-\alpha_b)J_N(t+d^*))dt,      \eqno (14)
$$
$$Q_b(x):=
:=\frac 1{\pi}\int_{0}^xq_b(t)\frac {\Pi_*(t)}{M}
(\alpha_bJ_N(t-d^*)+(1-\alpha_b)J_N(t+d^*))dt.\eqno (15)$$ В силу
лемм 5 и 2 $g_b$ и $Q_b$ суть непрерывные $2\pi$-периодические
функции. Лемма 5 немедленно влечет следующие соотношения
$$g_b(x)=0,\ \text {когда}\ x\in [-b, b];                       \eqno (16)$$
$$g_b\in\Bbb C^{(5)}\cap\bigtriangleup^{(1)}(Y);    \eqno (17)$$
$$\|g_b\|<1,\ \ \ \|g_b'\|<\tilde M;    \eqno (18)$$
$$
g_b'(x)=Q_b'(x),\quad \overline Mb<|x|\le\pi;
$$

$$\|g_b'-Q_b'\|\le
\tilde M\frac {(\overline M-1)b}2 \frac {(\overline M+1)b}2\frac
{\overline Mb}2\le \frac {\tilde M\overline M^3b^3}8;
$$
$$
\|g_b-Q_b\|\le \overline Mb\|g_b'-Q_b'\|\le \frac {\tilde
M\overline M^4b^4}8.\eqno (19)
$$
Легко видеть, что
$$
Q_b''(0)=-\frac 12\sin^2\frac b2 \frac{\Pi_*(0)}M\frac 1\pi
J_N(d^*) <-\frac 1{2\pi^2}\frac{m\tilde m}Mb^2,\eqno (20)
$$
$$
\| Q_b^{(k+1)}\|<M_k,
$$
где $M_k=$const, не зависит от $b$;
$$
\omega_k\left(g_b';t\right)\le 2^k\|g_b'-
Q_b'\|+t^k\|Q_b^{(k+1)}\| \le 2^{k-3}\overline M^3\tilde
Mb^3+t^kM_k. \eqno (21)
$$

\bigskip
{\bf П р и м е р 2.} {\it Для любых $k>3$ и $n\in\Bbb N$
существует функция  $f_2(x):=f_2(x;s,Y,n,k)$ такая, что
 $f_2\in\bigtriangleup^{(1)}(Y)\cap\Bbb C^{(1)}$ и
$$ E_n^{(1)}(f_2;Y)>B_Yn^{\frac k3 -1}\frac 1n\omega_k\left(f_2';\frac 1n
\right), \eqno (22) $$} где $B_Y=$const, зависит только от $Y$ и
$k$.

{\it Д о к а з а т е л ь с т в о.} Пусть $g(x):=g_b(x),\quad
Q(x):=Q_b(x),$ где $g_b(x)$ и $Q_b(x)$ определены соответственно
равенствами (14) и (15). Возьмем произвольный полином $\tau_n\in
\Bbb T_n\cap\bigtriangleup^{(1)}(Y)$, $n>s+2N-1$, и положим
$$R_n(x):=\tau_n(x)- Q(x).$$
Поскольку $\tau_n\in \bigtriangleup^{(1)}(Y)$, то $\tau_n''(0)\ge
0$, тогда в силу (20)
$$R_n''\left(0\right)=\tau_n''\left(0\right)- Q''\left(0
\right)\ge - Q''\left(0\right)\ge\frac {m\tilde
mb^2}{2\pi^2M}:=m_*b^2.$$ Применяя неравенство Бернштейна
$$\|R_n''\|\leq n^2\|R_n\|,
$$
получаем
$$m_*b^2\le |R_n''(0)|\le n^2\|R_n\|,$$
откуда с учетом (19),
$$m_*\frac {b^2}{n^2}\le\|R_n\|\leq
\|\tau_n-g\|+\|g-Q\|\le\|\tau_n-g\|+ \frac {\overline M^4\tilde
Mb^4}8\ ,$$ т.е.
$$\|\tau_n-g\|\ge m_*\frac {b^2}{n^2} -
\frac {\overline M^4\tilde Mb^4}8=m_*\frac{b^2}{n^2}\left( 1-\frac
{\overline M^4\tilde Mb^2n^2}{8m_*}\right).       \eqno (23) $$
Теперь обозначим
$$
b_n:=\left(\frac 1n\right)^{\frac k3},
$$
выберем $N_0$ из условий
$$
N_0>s+2N-1,\quad b_{N_0}<\frac 1{2\overline MN}, \quad \frac
{\overline M^4\tilde M}{8m_*}b^2_{N_0}N_0^2<\frac 12,
$$
и при всех $n\ge N_0$ положим
$$
f_2(x,s,Y,n,k):=g_{b_n}(x).
$$
Вследствие (17)
 $$f_2\in\bigtriangleup^{(1)}(Y)\cap\Bbb C^{(1)}.$$
Наконец, неравенство (22) следует из (21) и (23):
\begin{multline*} \frac{nE_n^{(1)}(f_2,Y)}{\omega_k(f_2';\frac
1n)}\ge \frac {n\|\tau_n-f_2\|}{\omega_k(f_2';\frac 1n)}\ge
m_*\frac {b_n^2}{n}\left( 1 -\frac {\overline M^4\tilde Mb_n^2n^2}
{8m_*}\right) :\\
 :\left( 2^{k-3}\overline M^3\tilde Mb_n^3 +\left(\frac 1n \right)^kM_k\right)\ge
\frac {m_*b_n^2}{2n}
:\left( 2^{k-3}\overline M^3\tilde Mb_n^3 +\left(\frac 1n \right)^kM_k\right)=:\\
=:B_Yn^{\frac k3 -1}. \end{multline*}

Для $n\ge N_0$ неравенство (22) доказано. Для случая $n<N_0$ оно
следует из неравенства $E_n^{(1)}(f_2;Y)\geq
E_{N_0}^{(1)}(f_2;Y)$. Пример доказан.

\smallskip

{\small \centerline{Л И Т Е Р А Т У Р А} \nopagebreak

1. {\it Jackson~D.} On Approximation by Trigonometric Sums and
Po\-ly\-no\-mi\-als// Trans. Amer. Math. Soc., {\bf 13} (1912),
S.~419--545.

2. {\it Zygmund~A.} Smooth Functions// Duke math. journ., {\bf 12}
(1945), 1, S.~46--76.

3. {\it Стечкин~С.\,Б.}  О наилучшем приближении периодических
функций  тригонометрическими полиномами// Доклады АН СССР,  {\bf
83} (1952), 5, C.~651--654.

4. {\it Копотун~К.\,А.} Равномерные оценки выпуклого приближения
многочленами// Мат. заметки, {\bf 51} (1992), 3, С.~35--46.

5. {\it Тиман~А.\,Ф.} Усиление теоремы Джексона о наилучшем
приближении непрерывных функций на конечном отрезке вещественной
оси//  ДАН СССР,  {\bf 78} (1951), 1, С.~17--20.

6. {\it Дзядык~В.\,К.} О приближении функций обыкновенными
многочленами на конечном отрезке вещественной оси// Изв. АН СССР,
сер. матем., {\bf 22} (1958), 3, С.~337--354.

7. {\it Freud~G.} Uber die Approximation Reelen Stetiger
Functionen Durch Gewohnliche Polinome// Math. Ann.,  {\bf 137}
(1959), 1, S.~17--25.

8. {\it Теляковский~С.\,А.} Две теоремы о приближении функций
 алгебраическими полиномами// Мат. сб., {\bf 79} (1966), 2, С.~252--265.

9. {\it Брудный~Ю.\,А.} Приближение функций алгебраическими
многочле нами// Изв. АН СССР, сер. матем., {\bf 32} (1968), 4,
С.~780--787.

10. {\it Lorentz~G.\,G., Zeller~K.\,L.}  Degree of Approximation
by Monotone Polynomials II// J. Approx. Theory,  {\bf 2} (1969),
3, S.~265--269.

11. {\it Шевчук~И.\,А.} Приближение многочленами и следы
непрерывных на отрезке функций. Киев:``Наукова думка'',  1992. 225
стр.

12. {\it Шведов~А.\,С.}  Комонотонное приближение функций
многочленами// ДАН СССР, {\bf 250} (1980), 1, 39--42.

13. {\it Дзядык~В.\,К.} Введение в теорию равномерного приближения
функций полиномами. M.:\newline "Наука",  1977. 512 стр.

\end{document}